\documentclass[letterpaper, 10 pt, conference]{ieeeconf}  %

\IEEEoverridecommandlockouts                              

\usepackage{amsmath} 
\usepackage{amssymb}  
\usepackage{tikz,pgfplots}
\usetikzlibrary{calc,intersections,shapes,arrows,positioning}
\pgfplotsset{compat=1.13}
\usepackage{graphicx}  
\usepackage{enumerate}
\usetikzlibrary{shapes.arrows}

\newtheorem{task}{Task}

\usepackage{caption}
\usepackage{subcaption}
\usepackage{siunitx}
\usepackage{color}
\definecolor{mycolor1}{rgb}{0.00000,0.5,0.75}

\title{\LARGE \bf
Learning References with Gaussian Processes in Model Predictive Control applied to Robot Assisted Surgery
}

\author{Janine Matschek$^{1}$, Tim Gonschorek$^{2}$, Magnus Hanses$^{3}$, Norbert Elkmann$^{3}$,  Frank Ortmeier$^{2}$, and Rolf Findeisen$^{1}$
\thanks{*This work was funded by the Federal Ministry of Education and Research in Germany within the Research Campus STIMULATE under Grant no. 13GW0095A. }
\thanks{$^{1}$Laboratory for Systems Theory and Automatic Control,
 	Otto-von-Guericke Univ. Magdeburg,
    {\tt\small \{janine.matschek, rolf.findeisen\}@ovgu.de }}%
\thanks{$^{2}$Chair of Software Engineering,
Otto-von-Guericke Univ. Magdeburg}%
\thanks{$^{3}$Fraunhofer-Institute for Factory Operation and Automation, Magdeburg}%
}

\begin{document}

\maketitle
\thispagestyle{empty}
\pagestyle{empty}

\begin{abstract}

One of the key benefits of model predictive control is the capability of controlling a system proactively in the sense of taking the future system evolution into account. However, often external disturbances or references are not a priori known, which renders the predictive controllers shortsighted or uninformed. Adaptive prediction models can be used to overcome this issue and provide predictions of these signals to the controller. In this work we propose to learn references via Gaussian processes for model predictive controllers. To illustrate the approach, we consider robot assisted surgery, where  a robotic manipulator needs to follow a learned reference position based on optical tracking measurements.

\end{abstract}

\section{INTRODUCTION}

Model predictive control (MPC) is an optimization based control strategy that can be applied to a variety of tasks including the control of constrained, nonlinear, multi-input multi-output systems. 
The involved optimization problem is defined for a prediction horizon, i.e.\ MPC is predicting the future evolution of system states based on a dynamical model of the system and evaluates its performance via a cost functional to be minimized over that horizon \cite{mayne2006robust}.
One of its mayor benefits is this predictive capability which makes MPC ``foresightful'': It can react proactively to known  upcoming events like changing references or external influences such as disturbances.
Different approaches to track time dependent signals with MPC exist
\cite{aydiner2016periodic, kern2009receding,limon2015tracking,limon2016mpc, maeder2010offset}.
Anyhow, the predictive capability can only be fully exploited if the future reference or disturbance is known a priori.  
This information is, however, not always a priori available.
Examples are changes of the wind conditions in UAV control, changing solar radiation and outside temperature in climate building control, or human drivers in a mixed autonomous car scenario.
If these external signals are a priori unknown but can at least be measured or observed and if insights from former conditions exist, a prediction model for those signals could be used.
For example wind and weather forecasts or the expected behaviour of human drivers could be modelled based on data from meteorological stations or traffic monitoring.
Using prediction models not only for the controlled dynamic system but also for these external signals allows for improving the predictive capabilities and consequently the expected performance of MPC. 

In this work Gaussian processes (\cite{Rasmussen2006}) are used to predict upcoming references online in a data-driven way.
Specifically, the Gaussian process (GP) will model the motion of an organ or goal structure in minimally invasive surgery based on optical tracking data.
Gaussian processes are machine learning techniques which allow to obtain a probabilistic system representation while handling noisy data and reduce overfitting \cite{Rasmussen2006}. 
It was shown in \cite{Klenske2016} that GPs can be used to predict quasi periodic signals by choosing specific covariance and mean functions.
These ideas are used to model an external reference signal which is provided to the model predictive controller.
We will follow this idea and extend it to a tracking MPC formulation similar to \cite{faulwasser11b} applied to medical robotics and motion compensation. 
Application examples of tracking MPC for motion compensation (without using GPs) in robot assisted surgery are given in \cite{gangloff2006model, arenbeck2017control, paluszczyszyn2013model}.

The main contribution of this work is the combination of learning based prediction models obtained via Gaussian processes with model predictive control applied to the field of medical robotics. More specific, the future evolution of an external measurement signal is learned online to be provided as a reference for the predictive controller. Both the predicted mean of the Gaussian process as well as its derivative are used to design the reference signal. As a specific application, robot assisted minimally invasive surgery is considered. The robot is controlled to precisely position a surgical instrument while following motions of the goal structures which originate e.g.\ from patient breathing. Hereby, the motion prediction should be smooth and result in prediction errors below the data noise while the controller needs to achieve a tracking position accuracy in submillimeter range and an orientational error below 0.1$^\circ$ to prevent from any harm of surrounding tissue and enable successful treatment.

The remainder of this paper is structured as follows:
Section~2 describes the problem setup including the Gaussian process and model predictive control formulations.
In Section~3 the control of a 7 DOF robotic manipulator for surgical applications is illustrating the proposed approach.
Section~4 summarises the stated results and gives possible directions for future work.

\section{PROBLEM DESCRIPTION}

Starting from a theoretical and a practical task description, we outline the basic concepts involved which are Gaussian processes, their use as reference generators as well as the resulting predictive control formulation.

\subsection{Objective}
We consider a dynamic system of the form
\begin{equation}
\dot x=f(x,u) \label{eq:system}
\end{equation}
where $x \in  \mathbb{R}^{n_x}$ and $u \in \mathbb{R}^{n_u}$ describe the states and inputs respectively.
Our goal is to control this system under the presence of state constraints $\mathcal X$ and input constraints $\mathcal U$ to follow a time varying reference $x_r(t): \mathbb{R} \mapsto \mathcal X$, i.e. $x$ should follow $x_r$ as good as possible.
We want to use model predictive control, as it allows to take  constraints directly into account and to find optimal inputs to the system minimizing the tracking error.
As the reference $x_r$ changes over time, we are interested in its future predictions to provide it to the controller.
To be able to predict the future reference evolution, we rely on data based machine learning instead of a priori known first principle models. For this purpose, reference generation involving Gaussian process regression  are derived predicting the reference  based on possibly noisy data.
Figure~\ref{fig:basic_idea} illustrates the basic setup.
The task can be summarized as follows:

\begin{task}
Consider system \eqref{eq:system} and given some data $\mathcal D$ obtain a reference generator and a controller such that the following holds
\begin{enumerate}[(i)]
\item The reference $x_r(t)$ is reachable, i.e. $x_r(t)\in \mathcal X$ and $\exists  r_u \in \mathcal U $ such that $\dot x_r(t)=f(x_r(t),r_u(t))$, and it is available for the controller during its prediction horizon.
\item The control error converges to zero, i.e. $\lim \limits_{t\rightarrow \infty} e(t)=\lim \limits_{t\rightarrow \infty}x(t)-x_r(t)=0$.
\item The state and input constraints are satisfied, i.e. $x\in \mathcal X$ and $u\in \mathcal U$.
\end{enumerate} 
\label{task}
\end{task}

\tikzset{
block/.style = {draw, fill=white, rectangle, minimum height=1.5cm, minimum width=1.3cm, align=center},
tmp/.style  = {coordinate}, 
sum/.style= {draw, fill=white, circle, minimum size=0.5cm},
input/.style = {coordinate},
output/.style= {coordinate},
pinstyle/.style = {pin edge={to-,thin,black},
>=latex
}
}

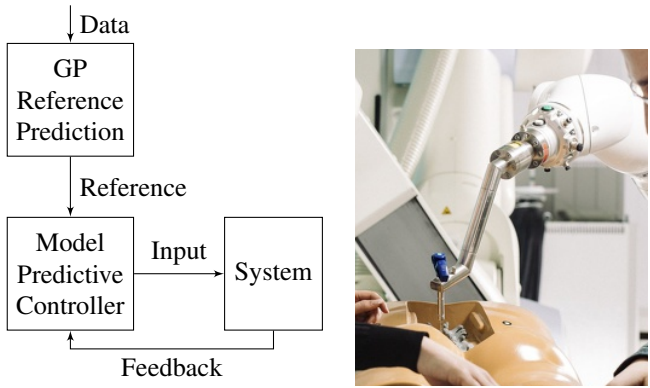
\begin{figure}[!tb]
\centering
\begin{subfigure}[t]{0.45\columnwidth}
\begin{tikzpicture}[auto, node distance=2.7cm,>=latex']
    \node [input] (input){};
    \node [block, below =0.5cm of input] (GP)
    		{GP \\ Reference\\ Prediction };
    \node [block, below =0.8cm of GP] (MPC) 
    		{Model\\ Predictive \\ Controller};
    \node [block, right of=MPC] (system) {System};
    
    \draw [->] (input) -- node{Data} (GP);
    \draw [->] (GP) -- node{Reference} (MPC);
    \draw [->] (MPC) -- node{Input} (system);
    \draw [->] (system) -- ++(0,-1) -| node[pos=0.25]
    				{Feedback} (MPC);      
\end{tikzpicture}
\caption{Basic idea of data-based reference prediction via GPs for model predictive control.}
\label{fig:basic_idea}
\end{subfigure} \hspace{0.5cm}
\begin{subfigure}[t]{0.45\columnwidth}
\begin{center}
\includegraphics[width=1\textwidth, trim=350 150 200 100,clip]{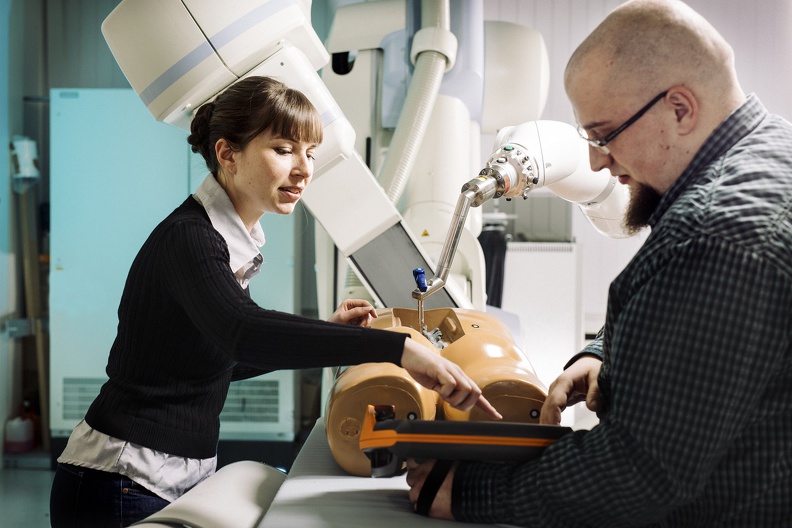}
\caption{Robot supported intervention on a phantom.}
\label{fig:robot}
\end{center}
\end{subfigure}
\caption{Theoretical and practical task illustration.}
\end{figure}

\subsection{Medical Task Description}
As application example we focus on the control of a seven degrees of freedom medical lightweight robot \cite{KUKA} that assists physicians during minimally invasive surgery, see also Figure~\ref{fig:robot}. More specific, our application test case is robot supported radio frequency ablation of tumors in the spine \cite{matschek2017mathematical, hanses2016robotic}. 
Hereby, needle shaped electrodes must be placed precisely inside the target area, which is in direct vicinity to the spinal cord.
Understandably enough, high position accuracy is needed which is defined as submillimeter Cartesian position accuracy and orientational errors below 0.1$^\circ$.
We propose to use a robotic manipulator as an assistance to the interventionalists which is controlled by a predictive controller to minimize positioning errors while taking constraints, as e.g.\ joint limitations or areas that should not be touched, into account.

As radio frequency ablation in the spine is usually performed in prone position (dorsal side up) motions of the patient, e.g.\ due to breathing, affect the Cartesian position of the goal structures. 
This motion can be tracked by optical tracking systems. However, control solely based on those tracking measurements  suffer from problems related to line-of-sight such as data loss and can only \emph{react} on the tracked motions. To overcome this, we propose to use this data to obtain a mathematical description, i.e.\ model, of the motion. This facilitates to compensate for data loss and enables a prediction of the likely future motion to allow the controller to be \emph{proactive}. 
The motion model should on one side be based on data (and e.g.\ handle noise in the measurements) as well as showing good extrapolation qualities as we want to predict into the future. 
Gaussian processes supply those characteristics if the mean and covariance functions are chosen to represent intrinsic system structure while predictions are performed by taking upcoming measurements into account.

\subsection{Gaussian Processes}

Gaussian Processes provide schemes for non-parametric, stochastic system identification. They have gained significant attention in the control community \cite{Rasmussen2006, Ostafew2016, kocijan2004gaussian, berkenkamp2015safe}.

Some reasons for this popularity are the limited amount of design decisions, their capability of dealing with noisy data, and the confidence interval that is provided by the GP which allows to investigate the quality of the obtained model.

A Gaussian distribution can be uniquely defined via a mean function $m(z): \mathbb{R}^{n_z} \mapsto \mathbb{R} $ and a symmetric, positive semi definite covariance function $k(z,z'): \mathbb{R}^{n_z}\times \mathbb{R}^{n_z} \mapsto \mathbb{R}_0^+$. It is written as 
\begin{equation*}
r(z) \sim \mathcal{GP} (m(z),k(z,z')).
\end{equation*} 
Here, $z, z' \in \mathbb{R}^{n_z}$ are the independent variables or inputs to the GP and the values of the process $r(z)$ at specific locations $z$ possess a normal distribution.
Via the covariance function (also called kernel) a GP relates similarities between the input variables to the similarity between the output variables. 
To make predictions  at test points a number of observations, a so called training data set $\mathcal{D}=\{z_i,r_i| i=1,2,\ldots,n_\mathcal{D}\}$, must be available. 
The posteriori distribution (prediction) for the test point $z_*$ is given by the prior belief (selected covariance and mean function) conditioned on the available data $\mathcal D$. The joint distribution of the training data output $r$ and the test data output $r_*$ (assuming $n_*$ test data point) can be expressed as
\begin{equation*}
\begin{pmatrix}
r\\r_*
\end{pmatrix} \sim \mathcal{GP}\left(   \begin{pmatrix}
m(z)\\m(z_*)
\end{pmatrix}  , \begin{pmatrix}
K(z,z)+\sigma^2_n I & K(z,z_*)\\ K(z_*,z) & K(z_*,z_*)
\end{pmatrix}\right).
\end{equation*}
Here, $\sigma^2_n$ is the variance of the measurement noise. The entries of the covariance matrix $K$ are calculated using the covariance function $k$. Specifically, $K(z,z)$ is of dimension $n_\mathcal{D}\times n_\mathcal{D}$ and specifies the covariance between all of the training data points, while $ K(z,z_*)$ and  $ K(z_*,z)$ (with dimensions $n_\mathcal{D}\times n_*$ and $n_* \times n_\mathcal{D}$, respectively) define the cross correlation between test and training data points. $ K(z_*,z_*)$ with dimension $n_* \times n_*$ is the auto covariance of the test data. 
Given this joint probability distribution, the conditional posterior distribution can be calculated via the posterior mean 
\begin{equation*}
m^+(r_*)= m(z_*)+K(z_*,z)(K(z,z)+\sigma^2_n I)^{-1}(r-m(z))
\end{equation*}
 and posterior covariance
 \begin{align*}
 cov^+(r_*)&=K(z_*,z_*)\\
 &-K(z_*,z)(K(z,z)+\sigma_n^2 I)^{-1}K(z,z_*).
 \end{align*}
 
We also use the derivative of the posterior with respect to the test points (the derivative of a GP is a GP \cite{kocijan2016modelling}) inside the reference generator. To this end, we rely only on the observations $z,r$ and not on their derivatives. 
Clearly, the covariance and mean functions contain a number of parameters, so called \emph{hyperparameters}. 
We calculate them via the maximization of the marginal logarithmic  likelihood. This determination of the hyperparameters refers to the actual \emph{learning} or \emph{training} of the GP.

\tikzset{
block/.style = {draw, fill=white, rectangle, minimum height=1.5cm, minimum width=3em, align=center},
tmp/.style  = {coordinate}, 
sum/.style= {draw, fill=white, circle, minimum size=0.5cm},
input/.style = {coordinate},
output/.style= {coordinate},
pinstyle/.style = {pin edge={to-,thin,black},
>=latex
}
}

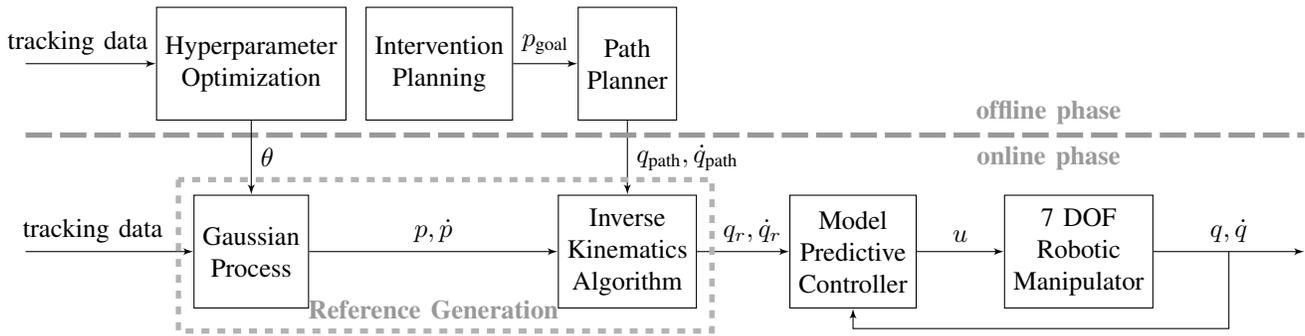
\begin{figure*}[!htb]
\centering
\begin{tikzpicture}[auto, node distance=3cm,>=latex']

    \node [input] (input){};
    \node [block, right of=input] (GP){Gaussian\\ Process};
    \node [block, right of=GP, node distance=5cm] (DIK) { Inverse\\ Kinematics\\ Algorithm};
    \node [block, right of=DIK] (MPC){Model\\ Predictive\\ Controller};
    \node [block, right of=MPC] (robot){7 DOF\\ Robotic\\ Manipulator};
    \node [output, right of=robot] (output){};
    
    \node [input, above of=input, node distance=1.55cm] (line_in){};
    \node [output, above of=output, node distance=1.55cm] (line_out){};
    
    \node [block, above of=DIK, node distance=2.5cm] (planner){Path\\ Planner};
    \node [block, above of=GP, node distance=2.5cm] (hyper){Hyperparameter\\ Optimization};
    \node [input, above of=input, node distance=2.5cm] (hyp_in){};
    \node [block, left of=planner, node distance=2.5cm] (inter){Intervention\\ Planning};

    \draw [->] (input) -- node[pos=0.4]{tracking data} (GP);
    \draw [->] (GP) --node[name=p]{$p,\dot{p}$} (DIK);
    \draw [->] (DIK) -- node[pos=0.6]{$q_r, \dot q_r$} (MPC);
        \draw [->] (planner) -- node{$q_\text{path},\dot q_\text{path}$}  (DIK);
    
    \draw [->] (MPC) --node{$u$} (robot); 
    \draw [->] (robot) -- node[name=y]{$q,\dot q$}(output);     
    \node [tmp, below of=y, node distance=1.3cm] (tmp1){};
    \draw [->] (y) |- (tmp1) -| (MPC);

    \draw [black!40!white,line width=0.7mm,dashed, dash pattern=on 10pt off 3pt] (line_in) --node[pos=0.8, below] {\textbf{online phase}}node[pos=0.8] {\textbf{offline phase}} (line_out);
    
    \draw [->] (hyp_in) -- node[pos=0.4]{tracking data} (hyper);
    \draw [->] (hyper) -- node{$\theta$} (GP);
    \draw [->] (inter) -- node{$p_\text{goal}$} (planner);
    
    \draw[black!30!white,line width=0.7mm,dashed] ($(GP.north west)+(-0.2,0.2)$)  rectangle ($(DIK.south east)+(0.2,-0.3)$);
    \node [ below of=p, node distance=1.05cm] (tmp1){\color{black!40!white}\textbf{Reference Generation}};

\end{tikzpicture}
\caption{System setup including the predictive controller which uses a reference prediction provided by the GP.}
\label{fig:loop}
\end{figure*}

\subsection{Reference Governors and Reference Generation}
In its most general form, a reference governor is used to provide a (modified) reference to a control loop. 
In \cite{gilbert1995discrete, bemporad1997nonlinear} nonlinear filtering techniques are used to design trackable references, i.e.\ references which satisfy the constraints.
Besides, model predictive controller based reference governors  \cite{klauvco2013mpc} as well as concepts of incorporating reference governors inside an MPC controller \cite{limon2015tracking,limon2016mpc} exist.
In the given setup the reference generation provides a reference to the predictive controller that fulfills the following conditions:
\begin{enumerate}
\item The reference should be based on data without overfitting it.
\item The reference should be known to the controller for the whole prediction horizon, i.e.  at $t_k, x_r(\tau)$ is known $\forall \tau \in [t_k,t_k+T_F]$. 
\item The reference should be reachable, i.e. it should fulfill the state constraints $x_r(t)\in \mathcal X$ and it should be followable given the system dynamics $\exists  r_u \in \mathcal{U} $ such that $\dot x_r(t)=f(x_r(t),r_u(t))$. 
\end{enumerate}
The first and second issue are solved by using Gaussian processes (and some suitable post processing mapping from $m^+(r_*)$ to $x_r,u_r$) which provide future predictions of the reference signal based on past observations.
As a first step to address the third requirement we assume that the original reference (a priori unknown, provided in terms of data) is fulfilling condition 3) and the GP is trained to map this underlying unknown function. Certainly, further considerations need to be taken into account to  guarantee condition 3) in the general case. However, for the considered example in Section~3 we are able to do so (the system is an integrator chain).

\subsection{Model Predictive Control}
\label{sec:mpc}

In general, MPC minimizes a given cost criterion over a (finite) prediction horizon to obtain an optimal input to the system under control. Once an optimal input signal for the considered prediction horizon is obtained, a part of it is applied to the real system and the optimization is repeated with a shifted horizon taking recent measurements into account.
Inside the optimization a model of the controlled system is used to predict the upcoming system evolution. Besides, MPC naturally handles restrictions on states and inputs by including them as constraints in the optimization procedure.
In \emph{tracking MPC} a time dependent reference should be tracked \cite{matschek2019nonlinear}. 
We assume that the reference generator provides this reference in terms of state and input references, $x_r(t)$ and $u_r(t)$, respectively.
A general formulation of tracking MPC is given by
\begin{subequations}
\begin{align}
	\min_{u}& {\int^{t+T_F}_{t} L\left(x(\tau), x_r(\tau), u(\tau),u_r(\tau)\right) d\tau}\\
	& + E(x(t+T_F),x_r(t+T_F))
	\intertext{subject to ( $\forall	\tau \in [t,t+T_F]$)}
	\dot{x}(\tau)&=f(x(\tau),u(\tau)), \quad x(t)=x_t\label{eq:MPC_dotx}, \\ 
	x(\tau)&\in \mathcal{X} \quad  
	u(\tau)\in \mathcal{U} \label{eq:MPC_conu}\\
	x(t+T_F)&\in \mathcal{E} \label{eq:MPC_conE}.
\end{align}
	\label{eq:ocp}
\end{subequations}
The cost functional to be minimized consists of the stage cost $L(\cdot,\cdot,\cdot,\cdot):\mathbb{R}^{n_x}\times\mathbb{R}^{n_x}\times\mathbb{R}^{n_u}\times\mathbb{R}^{n_u} \mapsto \mathbb{R}_0^+$ and the terminal cost $E(\cdot,\cdot):\mathbb{R}^{n_x}\times\mathbb{R}^{n_x} \mapsto \mathbb{R}_0^+$. 
The prediction horizon of the MPC problem is given by $T_F$, where due to the time dependency of the reference, the considered time span runs from the actual time $t$ to $t+T_F$ (even for time-invariant system dynamics).
The system dynamics are given by \eqref{eq:MPC_dotx} where $x\in \mathbb{R}^{n_x}$ are the states of the system, $u \in \mathbb{R}^{n_u}$ are the control inputs and $x_t \in \mathbb{R}^{n_x}$ are initial conditions.
Constraints on the states and inputs are expressed by \eqref{eq:MPC_conu}. 
A terminal constraint is given by \eqref{eq:MPC_conE}. 
The stage and terminal cost depend (besides of the states $x$ and inputs $u$) additionally on the state and input references $x_r\in \mathbb{R}^{n_x}$ and $u_r\in \mathbb{R}^{n_u}$. In our case, these values will be determined online via machine learning, i.e.\ based on a Gaussian process.  

\subsubsection*{Remark on Stability}
There exist several variants to guarantee asymptotic stability of model predictive control for tracking time-varying references \cite{matschek2019nonlinear}. One of them is to consider the time-varying error dynamics. In this case, stability criteria for time-varying MPC can be applied, see e.g. \cite{rawlings2017model, faulwasser11b}. 
One of the basic assumptions in those approaches is continuity of the error system which relates to differentiability of the reference.
If appropriate kernels are chosen, the mean value of a Gaussian process is continuously differentiable.
In case of fixed hyperparameters and training data, we can therefore directly apply the stability conditions of tracking MPC to our case.
However, if online learning is considered via the update of the training data set, the reference changes between two MPC executions such that differentiability of the closed loop reference might not be given. In this case, robust techniques might be used together with the fact, that the change of the reference induced by the GP training is bounded, see also \cite{maiworm2018stability}.  However, a detailed stability analysis for this case is beyond the scope of this work.

Given all previous considerations, we have defined the necessary ingredients to track external signals with MPC which are provided by data-based online learning via GPs.   This allows us to tackle the specific application example as outlined below.

\begin{figure}
\begin{center}
\input{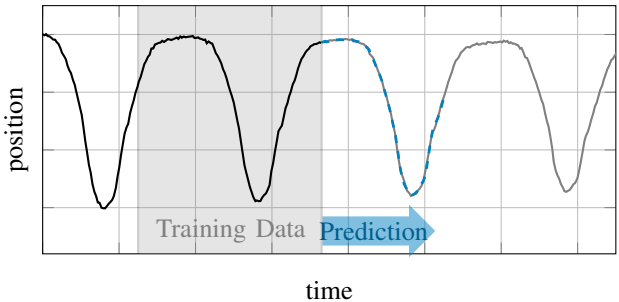}
\caption{Tracking data (black) is used to build a training data set of past observations (grey area) for GP prediction (dashed, blue) of future motion (grey line).}
\label{fig:tracking}
\end{center}
\end{figure}

\section{APPLICATION EXAMPLE: MEDICAL ROBOTICS}

This section describes the robotic application example and shows the achieved motion compensation via GPs and MPC.

\subsection{Specifically Proposed Controller Structure}

The overall setup of the proposed learning supported robot tracking control is shown in Figure~\ref{fig:loop}. 
In general there are two phases in this set up, namely an offline and an online phase.
During the \emph{offline phase} the intervention and path planning are performed. Based on a medical diagnosis the radiofrequency ablation is planned which results in the desired (Cartesian) electrode position and orientation $p_\text{goal}$.
Given this goal position and orientation of the electrode, the path planner determines a suitable needle insertion path in joint coordinates of the robot $q_\text{path},\dot q_\text{path}$. 
Additionally, the optimization of the hyperparameters $\theta$ for the Gaussian process is performed offline (based on patient specific interventional tracking data). This choice is justified by the clinical data, which showed no time-variance of the hyperparameters during an intervention, but variance for different patients.

Once the hyperparameters for the Gaussian process $\theta$ and the planned configurations of the robot $q_\text{path}$ and their derivatives $\dot q_\text{path}$ are obtained, the \emph{online phase} is executed.
Based on the tracking data that measures the actual position of the goal structure the training data of the Gaussian process is updated online. The GP predicts the future motion $p$ in terms of Cartesian positions and orientations (Euler angles) relative to the predefined goal $p_\text{goal}$. Additionally the derivative of this relative motion $\dot p$ is predicted.
To compensate this motion with minimum deviation from the pre-planned robot path $q_\text{path}$, an inverse kinematic algorithm is used. 
This results in a combined angular position $q_r$ which incorporates the planned path of electrode insertion as well as the motion compensation. 
The combined angular position $q_r$ is supplied as a reference to the model predictive controller, which generates optimal inputs for the robotic manipulator to track this reference.

\subsection{Reference Prediction by Gaussian Process}

The Gaussian process is used to model the motion of the goal structure due to patient motion that is mainly influenced by breathing during the intervention. 
The tracking data $\mathcal{D}=\{z_i,r_i| i=1,2,\ldots,n_\mathcal{D}\}$ consists of the time $t_i$ as independent variable $z_i$ and the relative Cartesian position and Euler angles as dependent variables $r_i$.
We use a receding horizon of $n_\mathcal{D}=200$ past observations (similar to moving horizon estimation), which is chosen to cover at least one full breathing cycle, see also Figure~\ref{fig:tracking}. 
Inspired by \cite{Klenske2016} the covariance function is assumed to be compounded by the squared exponential covariance function and a periodic covariance function, to model the quasi periodic nature of the motions. 
The test data for the Gaussian process $z_*$ represents the prediction horizon $[t,t+T_f]$ evaluated at discrete time instances in accordance to the implemented sampling time for the predictive controller.
The predicted mean value $m^+(z_*)$ is modelling the relative goal position and orientation $p$ over this prediction horizon.
Additionally, the derivative of the mean $\frac{\partial m^+(z_*)}{\partial z_*}$ supplies the velocity information $\dot p$.
As the dimension of the goal structure position consists of three Cartesian positions and three Euler angles, in total six values must be predicted. For simplicity of implementation, we have chosen to train six uncoupled Gaussian processes separately instead of taking the full interaction between the outputs into account. 
The variance of the prediction which indicates the quality of the prediction model is not used for the inverse kinematics algorithm. Rather, it is utilized as an additional safety layer by providing these  information online to the physician.

\subsection{Inverse Kinematics Algorithm}
To map the predicted relative Cartesian position and Euler angles to the joint space of the robot the following inverse kinematic algorithm \cite{siciliano2010robotics} is used
\begin{equation*}
\dot q_r =J_A^\dagger(q_r)(\dot p +K e)+(I -J_A^\dagger(q_r) J_A(q_r))\dot q_\text{path}
\end{equation*}
in combination with an integration. 
Here, $J_A^\dagger$ is the pseudo inverse of the analytic Jacobian of the manipulator and $I$ is an identity matrix of size $7\times 7$. The term $e=p-\tilde{f}_\text{fk}(q_r)$ defines the error between the predicted Cartesian position and orientation $p$ and the one based on the constructed angles $q_r$ via forward kinematics $\tilde{f}_\text{fk}$. So to say, $e$ can be viewed as a control error  where $q_r$ is fed back and $K$ can be seen as a controller matrix that should be chosen such that the inverse kinematic algorithm is asymptotically stable. Besides this feedback term, two feed forward terms are included in the algorithm to track the change  $\dot p $ and to assign a desired elbow motion via $\dot q_\text{path}$. 

\subsection{Model Predictive Controller}

The states of the dynamic system model $x=\begin{pmatrix}
q & \dot q
\end{pmatrix}^\top,$ are the seven robot angles $q$ and their angular velocities $\dot q$,
 i.e.\ the system size is $n_x=14$.  
The system model used in the predictive controller is a linear integrator chain 
 \begin{equation*}
 \dot x=\begin{pmatrix}
 0 & I \\0& 0
 \end{pmatrix}x+\begin{pmatrix}
 0 \\ I
 \end{pmatrix}u
 \end{equation*}
which can be obtained by using a full dynamic model of the robot including inertia, Coriolis, gravity, and friction terms to perform an input-output linearization. We have chosen to use the  input-output linearized model for predictions inside the MPC instead of the full dynamic model due to two reasons: First, the full dynamic model in highly nonlinear and complex which demands high computational power to solve the optimal control problem in real time with small sampling rates. Second, the motions we perform (motion compensation and electrode placement) are comparably small and slow, where the full potential of the complex dynamic system is neither fully exploited nor necessary. 

The constraints on the system states (angles and angular velocities) where chosen to be in line with the robot specifications but where not active.
We used a quadratic stage and terminal cost function $L=e_x^\top Q e_x+e_u^\top R e_u$ and $ E= e_u^\top Q_F e_u$ with the weighting matrices
$Q=\text{diag}[
1\cdot 10^6, 1\cdot 10^6, 1\cdot 10^6, 1\cdot 10^6, 1\cdot 10^6, 1\cdot 10^6, 1\cdot 10^6, 2\cdot 10^{-3}, 2\cdot 10^{-3}, 2 \cdot 10^{-3}, 2 \cdot 10^{-3}, 2\cdot 10^{-3}, 2 \cdot 10^{-3},2 \cdot 10^{-3}]$,
 $R=\text{diag} [2 \cdot 10^{-3},2 \cdot 10^{-3},2 \cdot 10^{-3},2 \cdot 10^{-3},2 \cdot 10^{-3},2 \cdot 10^{-3},2 \cdot 10^{-3}]$, and
$Q_F=1 \cdot 10^{-3}\text{diag}([1, 1, 1, 1, 1, 1, 1, 0.1, 0.1, 0.1, 0.1, 0.1, 0.1, 0.1])$.
The prediction horizon is $T_F=30\cdot T_s$, with a sampling time of $T_s=5\,$ms.

\subsection{Tracking Performance}

\begin{figure}
\begin{center}
\input{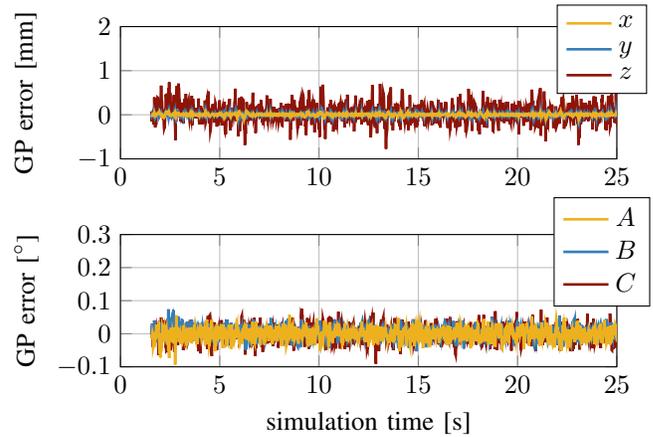}
\caption{Gaussian process prediction error.}
\label{fig:GP_error}
\end{center}
\end{figure}

\begin{figure}
\begin{center}
\input{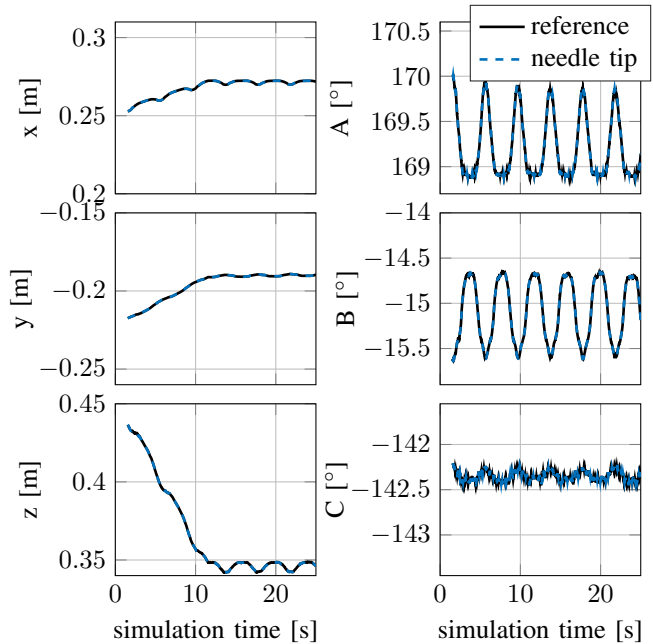}
\caption{Resulting robot motion following the composed reference path consisting of a 10cm Cartesian needle feed and the motion correction.}
\label{fig:Cart}
\end{center}
\end{figure}

\begin{figure}
\begin{center}
\input{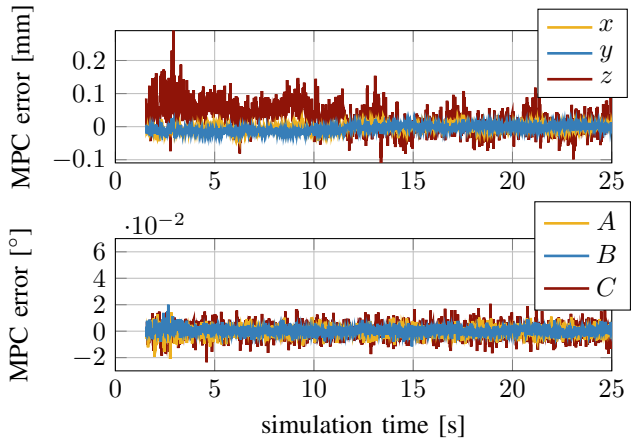}
\caption{Cartesian position and Euler angle errors.}
\label{fig:Cart_error}
\end{center}
\end{figure}

The closed loop (one step ahead) prediction quality of the GP is given in terms of the error in Cartesian position and orientation shown in Figure~\ref{fig:GP_error}, whose absolute values lie in sub-millimeter range and below $0.1^\circ$, respectively. This error, however, origins from the noise in the data, which we do not want to overfitt. The Gaussian process acts as a reference generator modelling the underlying true reference motion without overfitting the noisy data, which would not be trackable (and should not be tracked) by the robot.    
In Figure~\ref{fig:Cart} the composed Cartesian reference is depicted in black. It consists of the pre-planned entry path and the motion correction term predicted by the GP. The entry path in this case is a 10cm long motion along the needles start orientation, which leads to a dominant motion in Cartesian z-direction and smaller changes in Cartesian x- and y- direction. At the same time, the orientation along the path is fixed, such that the resulting orientational change is equal to the motion compensation plus the constant initial needle orientation. 
The model predictive controller is controlling the robot to follow this reference (defined in joint space via the  inverse kinematics algorithm) precisely, see Figure~\ref{fig:Cart} blue dashed line. The corresponding Cartesian position and orientational errors are depicted in Figure~\ref{fig:Cart_error}. Note, that these errors reflect the deviation of the needle tip to the reference provided by the GP. As the maximum errors lie below 0.3mm and 0.03$^\circ$, we have achieved the demanded precision.

%

\section{CONCLUSIONS}

In this work we have proposed to use Gaussian processes as reference generator for model predictive controllers. In cases where the reference is a priori unknown but can be measured or observed, a prediction model can provide useful information to the MPC to fully exploit its predictive capabilities and improve the control performance. We have demonstrated this approach with the control of a robotic manipulator applied in minimally invasive surgery. A needle placement task was performed while parts of the reference are not known a priori but need to be learned online.
Future work will focus on stability analysis for the use of GPs as reference governors as well as an experimental verification of the application example.


\bibliographystyle{IEEEtran}

\bibliography{ecc20_bibo}

\end{document}